\documentclass[reqno]{amsart}
\usepackage{amsmath}
\usepackage{amssymb}
\usepackage{mathabx}
\newcounter{makeconstant}
\newenvironment{makeconstant}%
{\refstepcounter{makeconstant}}%
{}
\def\mc#1{\begin{makeconstant}\label{#1}\end{makeconstant}}

\def\({\left(}
\def\){\right)}

\def\eul{{\rm e}}
\def\dee{\thinspace{\rm{d}}}

\DeclareMathOperator{\ord}{ord}
\DeclareMathOperator{\lcm}{lcm}

\let\divides\mid

\let\notdivides\nmid

\def\mid{:}

\newcommand{\fix}{\operatorname{\mathsf F}}

\newcommand{\orbit}{\operatorname{\mathsf O}}
\newcommand{\mertens}{\operatorname{\mathsf M}}

\theoremstyle{plain}
\newtheorem{theorem}{Theorem}
\newtheorem{lemma}[theorem]{Lemma}

\newtheorem{proposition}[theorem]{Proposition}%

\theoremstyle{definition}

\newtheorem{remark}[theorem]{Remark}
\newtheorem{remarks}[theorem]{Remarks}

\renewcommand{\le}{\leqslant}
\renewcommand{\ge}{\geqslant}
\renewcommand{\epsilon}{\varepsilon}

\def\bigo{\operatorname{O}}    
\def\littleo{\operatorname{o}} 

\def\NN{{\mathbb N}}
\def\PP{{\mathbb P}}

\def\RR{{\mathbb R}}

\def\Ll{{\mathcal L}}
\def\Nn{{\mathcal N}}
\def\Mm{{\mathcal M}}

\def\d{\delta}
\def\e{\varepsilon}
\def\th{\theta}

\def\lvert{\left\vert}
\def\rvert{\right\vert}

\def\ds{\displaystyle}

\begin{document}

\title{Automorphisms with exotic orbit growth}
\author{Stephan Baier}\address{(SB) Mathematisches Institut,
Georg-August Universit\"at, G\"ottingen, Germany}
\author{Sawian Jaidee}\address{(SJ) Department of Mathematics,
123 Mittraphab Road, Khon Kaen University 40002, Thailand}
\author{Shaun Stevens}
\author{Thomas Ward}\address{(SS \& TW) School of Mathematics,
University of East Anglia, Norwich NR4 7TJ, UK}
\date{\today}

\thanks{SS and TW particularly thank Tim Browning for his endless patience
in answering our Analytic Number Theory questions, and for putting us in touch with SB}
\thanks{SS supported by EPSRC grant EP/H00534X/1}

\begin{abstract}
We exhibit continua on two different growth scales in the
dynamical Mertens' theorem for ergodic automorphisms of
one-dimensional solenoids.
\end{abstract}

\maketitle

\section{Introduction}
\label{Suffering long time angels enraptured by Blake}

Automorphisms of compact metric groups provide a simple family
of dynamical systems with additional structure, rendering them
particularly amenable to detailed analysis. On the other hand,
they are rigid in the sense that they cannot be smoothly
perturbed, and for a fixed compact metric group the group of
automorphisms is itself countable and discrete. Thus it is not
clear which, if any, of their dynamical properties can vary
continuously. The most striking manifestation of this is that
it is not known if the set of possible topological entropies is
countable or is the set~$[0,\infty]$ (this question is
equivalent to Lehmer's problem in algebraic number theory; see
Lind~\cite{MR0346130} or the monograph~\cite{MR2000e:11087} for
the details). The possible exponential growth rates for the
number of closed orbits is easier to decide, and it is shown
in~\cite{MR2085157} that for any~$C\in[0,\infty]$ there is a
compact group automorphism~$T:X\to X$ with
\begin{equation}\label{thesunwassettingoveravalon}
\textstyle\frac1n\log\fix_T(n)\longrightarrow C
\end{equation}
as~$n\to\infty$, where~$\fix_T(n)=\lvert\{x\in X\mid
T^nx=x\}\rvert$. Unfortunately, the systems constructed to
achieve this continuum of different growth rates are
non-ergodic automorphisms of totally disconnected groups, and
so cannot be viewed as natural examples from the point of view
of dynamical systems. It is not clear if a result
like~\eqref{thesunwassettingoveravalon} is possible within the
more natural class of ergodic automorphisms on connected
groups, unless~$C$ is a logarithmic Mahler measure (in which
case there is a toral automorphism that achieves this).

Our purpose here is to indicate some of the diversity that is
nonetheless possible for ergodic automorphisms of connected
groups, for a measure of the growth in closed orbits that
involves more averaging than
does~\eqref{thesunwassettingoveravalon}. To describe this,
let~$T:X\to X$ be a continuous map on a compact metric space
with topological entropy~$h=h(T)$. A closed orbit~$\tau$ of
length~$\lvert\tau\rvert=n$ is a
set~$\{x,T(x),T^2(x),\dots,T^n(x)=x\}$ with cardinality~$n$.
Following the analogy between closed orbits and prime numbers
advanced by work of Parry and Pollicott~\cite{MR727704} and
Sharp~\cite{MR1139566}, asymptotics for the expression
\[
\mertens_T(N)=\sum_{\lvert\tau\rvert\le
N}\frac{1}{e^{h(T)\lvert\tau\rvert}}
\]
may be viewed as dynamical analogues of Mertens' theorem. The
expression~$\mertens_T(N)$ measures in a smoothed way the
extent to which the topological entropy reflects the
exponential growth in closed orbits or periodic points. A
simple illustration of how~$\mertens_T$ reflects this is to
note that\mc{expgrowthinfix} if
\[
\fix_T(n)=C_{\ref{expgrowthinfix}}\eul^{hn}+\bigo\left(\eul^{h'n}\right)
\]
for some~$h'<h$,
then\mc{constantinexpgrowthinfix}
\[
\mertens_T(N)=C_{\ref{expgrowthinfix}}\sum_{n=1}^{N}\frac1n+
C_{\ref{constantinexpgrowthinfix}}+\bigo\left(\eul^{-h''N}\right)
\]
for some~$h''>0$
(see~\cite{MR2486259}).

Writing~$\orbit_T(n)$ for the number of closed orbits of
length~$n$, we have
\[
\fix_T(n)=\sum_{d\divides n}d\orbit_T(d)
\]
and hence
\begin{equation}\label{Burn out the dross innocence captured again}
\orbit_T(n)=\frac{1}{n}\sum_{d\divides n}
\mu\left(\textstyle\frac{n}{d}\right)\fix_T(d)
\end{equation}
by M{\"o}bius inversion.

For continuous maps on compact metric spaces, it is clear that
all possible sequences arise for the count of orbits
(see~\cite{MR1873399}; a more subtle observation is that the
same holds in the setting of~$C^{\infty}$ diffeomorphisms of
the torus by a result of Windsor~\cite{MR2422026}). For
algebraic dynamical systems the situation is far more
constrained, and it is not clear how much freedom there is in
possible orbit-growth rates. Our purpose here is to exhibit two
different continua of growth rates, on two different speed
scales:
\begin{itemize}
\item for any~$\kappa\in(0,1)$ there is
an automorphism~$T$
of a one-dimensional compact metric group with~$\mertens_T(N)\sim
\kappa\log N$;
\item for any~$\delta\in(0,1)$ and~$k>0$ there is
an automorphism~$T$ of a one-dimensional compact metric group
with~$\mertens_T(N)\sim k(\log N)^{\delta}$.
\end{itemize}
While this plays no part in the argument, it is worth noting that
there is a complete divorce between the topological entropy
and the growth in closed orbits of these examples -- they all
have topological entropy~$\log2$.

\section{The systems studied}

We will study a family of endomorphisms (or automorphisms) of one-dimensional solenoids,
all built as isometric extensions of the circle-doubling map. To describe
these, let~$\mathbb P$ denote the set of rational primes, and associate to
any~$S\subset\mathbb P$ the ring
\[
R_S=\{r\in\mathbb Q\mid
\lvert r\rvert_p\le 1\mbox{ for all }p\in\mathbb P\setminus S\},
\]
where~$\lvert \cdot\rvert_p$ denotes the normalized~$p$-adic
absolute value on~$\mathbb Q$, so that~$|p|_p=p^{-1}$. Thus,
for example,~$R_{\emptyset}=\mathbb Z$,~$R_{\{2,3\}}=\mathbb
Z[\frac16]$, and~$R_{\mathbb P}=\mathbb Q$. Let~$T=T_S$ denote
the endomorphism of~$\widehat{R_S}$ dual to the map~$r\mapsto
2r$ on~$R_S$. This map may be thought of as an isometric
extension of the circle-doubling map, with topological
entropy~$h(T_S)=\sum_{p\in S\cup\{\infty\}}
\max\{\log\lvert2\rvert_p,0\}=\log 2$ (see~\cite{MR961739} for
an explanation of this formula, and for the simplest examples
of how the set~$S$ influences the number of periodic orbits).
Each element of~$S$ destroys some closed orbits, by lifting
them to non-closed orbits in the isometric extension;
see~\cite{MR2180241} for a detailed explanation in the
case~$S=\{2,3\}$. This is reflected in the formula for the
count of periodic points in the system,
\begin{equation}\label{Standing on the beach at sunset all the boats}
F_{T_S}(n)=(2^n-1)\prod_{p\in S}\lvert 2^n-1\rvert_p
\end{equation}
(see~\cite{MR1461206} for the general formula being used here),
showing that each inverted prime~$p$ in~$S$ in the dual group~$R_S$
removes the~$p$-part of~$(2^n-1)$ from the total count of all points
of period~$n$. The effect of each inverted prime in~$R_S$ on
the count of closed orbits via the
relation~\eqref{Burn out the dross innocence captured again}
is more involved.

We write~$\lvert x\rvert_S=\prod_{p\in S}\lvert x\rvert_p$ for
convenience, and since we will be using the same underlying map
throughout, we will replace~$T=T_S$ by the parameter~$S$
defining the system in all of the expressions from
Section~\ref{Suffering long time angels enraptured by
  Blake}. There are then three natural cases: the `finite' case
with~$\lvert S\rvert<\infty$ and the `co-finite' case with~$\lvert\mathbb
P\setminus S\rvert<\infty$, together producing countably
many examples, and the more complex remaining `infinite
and co-infinite' case. A special case of the results in~\cite{EMSW} is
that for~$S$ finite we have
\[
\mertens_{T_S}(N)=\mertens_S(N)=\sum_{\lvert\tau\rvert\le
N}\frac{1}{\eul^{h\lvert\tau\rvert}}= k_S\log N + C_S +\bigo\left(N^{-1}\right),
\]
for some~$k_S\in(0,1]\cap\mathbb Q$ and constant~$C_S$. For
example,~\cite[Ex.~1.5]{EMSW} shows
that
\[
k_{\{3,7\}}=\frac{269}{576}.
\]
Here we continue the analysis further, showing the following
theorem.

\begin{theorem}\label{thm:dense}
The set of possible values of the constant~$k_S$ with
\[
\mertens_S(N)=k_S\log N + C_S +\bigo\left(N^{-1}\right),
\]
as~$S$ varies among the finite subsets of~$\mathbb P$,
is dense in~$[0,1]$.
\end{theorem}

If~$S$ is infinite, then more
possibilities arise, but with less
control of the error terms.

\begin{theorem}\label{thm:onto} For any~$k\in(0,1)$, there is an
infinite co-infinite subset~$S$ of~$\mathbb P$ with
\[
\mertens_{S}(N)\sim k\log N.
\]
\end{theorem}
We also give explicit examples of sets~$S$ for which the value of~$k$
arising in Theorem~\ref{thm:onto} is transcendental.

The co-finite case is very different in that~$\mertens_S(N)$
converges as~$N\to\infty$; other orbit-counting asymptotics better
adapted to the polynomially bounded orbit-growth present in these systems
are studied in~\cite{EMSW2} and~\cite{Sawian}.
The following result is more surprising, in that a positive proportion
of primes may be omitted from~$S$ while still destroying so many orbits
that~$\mertens_S(N)$ is bounded.

\begin{proposition}\label{prop:zero}
There is a subset~$S$ of~$\mathbb P$ with natural density in~$(0,1)$
such that
\[
\mertens_{S}(N)=C_S+\bigo\left(N^{-1}\right).
\]
\end{proposition}
In fact there are such sets~$S$ with arbitrarily small non-zero
natural density.

While it seems hopeless to describe fully the range of possible growth
rates for~$\mertens_S(N)$ as~$S$ varies, we are able to exhibit many
examples whose growth lies strictly between that of the
examples in Theorem~\ref{thm:onto}
and that of the examples in Proposition~\ref{prop:zero}.

\begin{theorem}\label{thm:logdelta} For any~$\delta\in(0,1)$ and
any~$k>0$, there is a subset~$S$ of~$\mathbb P$ such that
\[
\mertens_{S}(N)\sim k\(\log N\)^\delta.
\]
\end{theorem}

We also find a family of examples whose growth lies between
that of the examples in Theorem~\ref{thm:logdelta} and of those
in Proposition~\ref{prop:zero}.

\begin{theorem}\label{thm:loglog} For any~$r\in\NN$ and any~$k>0$,
there is a subset~$S$ of~$\mathbb P$ with
\[
\mertens_{S}(N)\sim k\(\log\log N\)^r.
\]
\end{theorem}
Moreover, it is possible to achieve growth asymptotic to any
suitable function growing slower than~$\log\log N$. A byproduct
of the constructions for Theorems~\ref{thm:dense}
and~\ref{thm:loglog} gives sets~$S$ such that
both~$\mertens_{S}(N)$ and~$\mertens_{{\mathbb P\setminus
S}}(N)$ are~$\littleo\(\log N\)$.

The idea behind the proofs of all these result is rather similar.
We choose our set of primes~$S$ so that it
is easy to isolate a subseries of dominant
terms in~$M_{S}(N)$ in such a way that the sum of the remaining terms
converges, usually quickly (controlling this
rate governs the error terms). We describe a general framework
for dealing with such sets, and then the sets used to
carry out the constructions are
defined by arithmetical criteria relying on properties
of the set of primes~$p$ for which~$2$
has a given multiplicative order modulo~$p$.

\subsubsection*{Notation}

From~\eqref{Standing on the beach at sunset all the boats}, we
have~$\mertens_S(N)=\mertens_{S\cup\{2\}}(N)$, for any set~$S$: thus,
without loss of generality, we make the standing assumption
that~$2\notin S$. We will use various global
constants~$C_1,C_2,\dots$, each independent of~$N$
and numbered
consecutively. The symbols~$C$ and~$C_S$ denote local constants
specific to the statement being made at the time.
For an odd prime~$p$, denote by~$m_p$
the multiplicative order of~$2$ modulo~$p$; for a set~$T$ of
primes, we write~$m_T= \lcm\{m_p\mid p\in T\}$.
We will also use Landau's big-$\bigo$ and little-$\littleo$ notation.

\section{Asymptotic estimates}
\label{In a new world, me and you, girl}

By~\eqref{Burn out the dross innocence captured again}
and~\eqref{Standing on the beach at sunset all the boats} we have
\begin{align*}
\mertens_S(N)&=\sum_{n\le N}\frac{1}{n2^n}
\sum_{d\divides n}\mu\left(\textstyle\frac{n}{d}\right)
\lvert 2^d-1\rvert
\times
\lvert 2^d-1\rvert_S\\
&=\underbrace{\sum_{n\le N}\frac{\lvert 2^n-1\rvert_S}{n}}_{=F_S(N)}+
R_S(N),
\end{align*}
where the last equation defines both~$F_S(N)$ and~$R_S(N)$.

\begin{lemma}\label{All the boats keep moving slow}
$R_S(N)=C_S+\bigo\left(2^{-N/2}\right)$.
\end{lemma}

\begin{proof}
By definition,~$R_S(N)$ is the sum of two terms,
\[
R_S(N) = -\sum_{n=1}^{N}\frac{\lvert 2^n-1\rvert_S}{n2^n}\ +\
\sum_{n=1}^{N}\frac{1}{n2^n}\sum_{d\divides n,d<n}\mu\left(
\textstyle\frac{n}{d}\right)(2^d-1)\lvert 2^d-1\rvert_S.
\]
Since $\ds\frac{\lvert 2^n-1\rvert_S}{n2^n}\le \frac 1{2^n}$
and $\ds\frac{1}{n2^n}\sum_{d\divides n,d<n}(2^d-1)\lvert
2^d-1\rvert_S\le \frac 1{2^{n/2}}$, both sums converge
(absolutely) so $R_S(N)$ converges to some $C_S$. Moreover,
\[
\lvert
R_S(N)-C_S\rvert
\le
\sum_{n=N+1}^{\infty}\frac{1}{2^n}+\sum_{n=N+1}^{\infty}\frac{1}{2^{n/2}}
=\bigo\left(2^{-N/2}\right).
\]
\end{proof}

Thus we think of~$F_S(N)$ as a dominant term, and much of our effort will
be aimed at understanding how~$F_S(N)$ behaves as a function of~$S$,
which starts with understanding the arithmetic of~$2^n-1$. The main
tool here is the elementary observation that, for~$p$ a prime
and~$n\in\NN$,
\begin{equation}\label{eqn:ordp2n-1}
\ord_p\(2^n-1\) = \begin{cases}
\ord_p\(2^{m_p}-1\)+\ord_p(n) &\mbox{ if }m_p\divides n, \\
0 &\mbox{ otherwise,} \end{cases}
\end{equation}
where~$\ord_p(n)$ denotes the index of the highest power of~$p$
dividing~$n$, so that~$|n|_p=p^{-\ord_p(n)}$. In particular, if~$T$ is
a finite set of primes and~$n\in\mathbb N$, then
\begin{equation}\label{eqn:|2n-1|T}
\lvert 2^{nm_T}-1 \rvert_T = \lvert 2^{m_T}-1 \rvert_T \lvert n\rvert_T.
\end{equation}


In the proof of~\cite[Proposition~5.3]{EMSW}, a recipe is given
for computing the coefficient of~$\log N$ in the asymptotic
expansion of~$F_S(N)$, when~$S$ is finite. This is based on an
inclusion-exclusion argument, splitting up the sum~$F_S(N)$
according to the subsets of~$S$. The disadvantage of this
approach is that many subsets of~$S$ can lead to an empty sum:
in principle, the splitting works for infinite~$S$
(since~$F_S(N)$ is anyway a finite sum) but then the
decomposition of the sum falls into an uncountable number of
pieces. Here we take a different approach, splitting up the
expressions arising according to the values~$m_T$, for~$T$ a
finite subset of~$S$, rather than according to the subsets~$T$.
In this setting, several different subsets~$T$ may give the
same value for~$m_T$.

To this end, we set~$\Mm_S=\{m_p\mid p\in S\}$
and denote by~$\widebar\Mm_S$ its closure under taking least common multiples:
\[
\widebar\Mm_S=\{\lcm(\Mm')\mid \Mm'\subseteq\Mm_S\} =
\{m_T\mid T\subseteq S\}.
\]
For~$n\in\NN$, we put
\[
\bar m_n = \max\{\bar m\in\widebar\Mm_S\mid \bar m\divides n\} =
\lcm\{m_p\in\Mm_S \mid m_p\divides n\}.
\]
Then, for~$\bar m\in\widebar\Mm$, set
\[
\NN_{\bar m}=\left\{n\in\NN\mid\bar m_n=\bar
m\right\}
\]
and
\[
S_{\bar m}=\left\{p\in S \mid m_p\divides \bar m\right\}.
\]
Note that the sets~$S_{\bar m}$ are finite, even if~$S$ is
not. Then, using~\eqref{eqn:|2n-1|T}, we get
\begin{eqnarray}
F_{S}(N) &=& \sum_{\bar m\in\widebar\Mm_S}
\sum_{\genfrac{}{}{0pt}{}{n\le N}{n\in\NN_{\bar m}}}
  \frac{\lvert2^n-1\rvert_{S}}n \notag \\
&=& \sum_{\bar m\in\widebar\Mm_S} \frac{\lvert 2^{\bar m}-1\rvert_{S_{\bar m}}}{\bar m}
\sum_{\genfrac{}{}{0pt}{}{n\le N/\bar m}{m_p\notdivides n\bar m\text{ for }p\in
    S\setminus S_{\bar m}}}
  \frac{\lvert n\rvert_{S_{\bar m}}}n.
\label{eqn:howhorrible0}
\end{eqnarray}
The asymptotic behaviour of these inner sums can be computed, at least
in principle, using the results of~\cite[\S5]{EMSW}. However, for a
general set of primes~$S$, this is
cumbersome, so we
will specialize to sets which are easier to deal with. In this, we are
motivated by the next lemma, which follows
one of the many paths used to prove
Zsigmondy's theorem (see~\cite[{\S~8.3.1}]{MR2135478} for the details).

\begin{lemma}\label{Close your eyes and you'll find}
Fix~$n\in\mathbb N$, and let~$S\subset\PP$ be a set of primes containing
$\{p\in\PP\mid m_p=n\}$. Then
\[
\lvert2^n-1\rvert_S\le\frac{n}{2^{\phi(n)-2}}.
\]
\end{lemma}

Recall that a divisor of~$2^n-1$ is~\emph{primitive} (in the
sequence~$\(2^n-1\)_{n\ge1}$) if it has no common factor
with~$2^m-1$, for any~$m$ with~$1\le m<n$. Thus~$\{p\in\PP\mid
m_p=n\}$ is the set of primitive prime divisors of~$2^n-1$.
This set is finite, since~$m_p\ge\log_2p$, but may be large --
for example~$m_{233}=m_{1103}=m_{2089}=29$.
Schinzel~\cite{MR0143728} proved that there are infinitely
many~$n$ for which this set contains at least~$2$ elements, but
it seems that not much more is known about it in general.

\begin{proof}[Proof of Lemma~\ref{Close your eyes and you'll find}]
Writing~$(2^n-1)^*$ for the maximal primitive divisor of~$2^n-1$, we
certainly have
\begin{equation}\label{When the rhythm takes over your mind}
\lvert 2^n-1\rvert_S^{-1}\ge(2^n-1)^*.
\end{equation}
By factorizing~$x^n-1$ we have
\begin{equation}\label{eqn:2nphi}
2^n-1=\prod_{d\divides n}\Phi_d(2),
\end{equation}
where~$\Phi_d$ is the~$d$th cyclotomic polynomial. It follows that~$(2^n-1)^*$
is a factor of~$\Phi_n(2)$. If a prime~$p$
divides~$\gcd\(\Phi_n(2),\Phi_d(2)\)$ for some~$d\divides n$
with~$d<n$, then~$p\divides 2^d-1$. Then, from~\eqref{eqn:ordp2n-1},
\[
\ord_p(2^n-1)=\ord_p(2^d-1)+\ord_p(n/d)
\]
and, from~\eqref{eqn:2nphi},
\[
\ord_p(2^n-1)\ge\ord_p(2^d-1)+\ord_p(\Phi_n(2))\ge\ord_p(2^d-1)+1,
\]
so in particular~$p$ divides~$n/d$; therefore~$p$ divides~$n$
and~$d$ divides~$n/p$, so~$p$ divides
\[
(2^{n/p}-1).
\]
Moreover
\[
\ord_p(2^n-1)=\ord_p(2^{n/p}-1)+1
\]
and
\[
\ord_p(2^n-1)\ge\ord_p(2^{n/p}-1)+\ord_p(\Phi_n(2)),
\]
so in fact~$\ord_p(\Phi_n(2))=1$. Thus
$\gcd\left(\Phi_n(2),\prod_{d\divides n,d<n}\Phi_d(2)\right)$ divides
$\prod_{p\divides n} p$, which is at most $n$, and
\begin{equation}\label{We'll be swaying, music playing}
(2^n-1)^*\ge\Phi_n(2)/n.
\end{equation}
On the other hand, by M{\"o}bius inversion applied to~\eqref{eqn:2nphi},
\[
\Phi_n(2)=\prod_{d\divides n}(2^d-1)^{\mu(n/d)}
\]
so
\[
\log(\Phi_n(2))
=
\phi(n)\log(2)+\sum_{d\divides n}\mu(n/d)\log(1-2^{-d}),
\]
where~$\phi(n)$ is the Euler totient function. Now, using the Taylor
expansion for the logarithm,
\[
\lvert\sum_{d\divides n}\mu(n/d)\log(1-2^{-d})\rvert
\le \sum_{d\divides n}\sum_{j=1}^{\infty}\frac{2^{-jd}}{j}
= \sum_{j=1}^{\infty}\frac{2^{-j}}{j}\sum_{d\divides n}2^{-j(d-1)}
\le 2\log 2,
\]
so~$\Phi_n(2)\ge 2^{\phi(n)-2}$ and the result follows
by~\eqref{When the rhythm takes over your mind}
and~\eqref{We'll be swaying, music playing}.
\end{proof}

This lemma will be used as follows. Instead of starting with a set~$S$
of primes, we begin with~$\Mm$ a subset of $\mathbb N$ and put
\[
S_\Mm=\{p\in\PP\mid m_p\in \Mm\}.
\]
Then
\[
\sum_{n\le N}\frac{\lvert2^n-1\rvert_{S_\Mm}}n=
\underbrace{\sum_{\genfrac{}{}{0pt}{}{n\le N}{n\notin\Mm}}
\frac{\lvert2^n-1\rvert_{S_\Mm}}n}_{=D_{S_\Mm}(N)}+
\underbrace{\sum_{\genfrac{}{}{0pt}{}{n\le N}{n\in\Mm}}
\frac{\lvert2^n-1\rvert_{S_\Mm}}n}_{=Q_{S_\Mm}(N)},
\]
and, by Lemma~\ref{Close your eyes and you'll find},\mc{newforMcase}
\[
Q_{S_\Mm}(N)\le
C_{\ref{newforMcase}}\sum_{n\le N}\frac{1}{2^{\phi(n)}},
\]
which converges since $\phi(n)\ge\sqrt{n}$, for $n\ge 6$. Moreover,
the same observation shows that\mc{newforMcase2}
\[
Q_{S_\Mm}(N) = C_{\ref{newforMcase2}}+\bigo\left(2^{-\sqrt{N}}\right).
\]
Thus the asymptotic behaviour is governed by the dominant
term~$D_{S_\Mm}(N)$.
From Lemma~\ref{All the boats keep moving slow}, we
get\mc{newforMcase3}
\begin{equation}\label{eqn:mertensM}
\mertens_{S_\Mm}(N)=\sum_{\genfrac{}{}{0pt}{}{n\le N}{n\notin\Mm}}
\frac{\lvert2^n-1\rvert_{S_\Mm}}{n}+C_{\ref{newforMcase3}}
+\bigo\left(2^{-\sqrt{N}}\right).
\end{equation}
All our examples will take this form.


\begin{remarks}
(i)
We have set up maps~$S\mapsto \Mm_S$ and~$\Mm\mapsto S_\Mm$ between
the power sets of~$\PP\setminus\{2\}$ and~$\NN$, which are
order-preserving for inclusion. It is easy to check
that~$S_{\Mm_S}\supseteq S$, while~$S_{\Mm_{S_\Mm}}=S_\Mm$. Similarly,
we have~$\Mm_{S_\Mm}=\Mm\setminus\{1,6\}$, since all but the first and
sixth terms of the Mersenne sequence have primitive divisors,
and~$\Mm_{S_{\Mm_S}}=\Mm_S$. In particular, we can apply the
decomposition~\eqref{eqn:howhorrible0} of~$\mertens_{S_\Mm}(N)$ in
tandem with~\eqref{eqn:mertensM}. When we do so, we will
replace~$\widebar\Mm_{S_\Mm}$ by the closure~$\widebar\Mm$ of~$\Mm$
under least common multiples to get\mc{newforMcase5}
\begin{equation}
\mertens_{S_\Mm}(N) =
\sum_{\bar m\in\widebar\Mm} \frac{\lvert 2^{\bar m}-1\rvert_{S_{\bar m}}}{\bar m}
\sum_{\genfrac{}{}{0pt}{}{n\le N/\bar m,\ n\bar m\notin\Mm}
{m_p\notdivides n\bar m\text{ for }p\in S_\Mm\setminus S_{\bar m}}}
  \frac{\lvert n\rvert_{S_{\bar m}}}n
+ C_{\ref{newforMcase5}} + \bigo\(2^{-\sqrt N}\).
\label{eqn:howhorrible}
\end{equation}
(ii)
For a general set~$S$, let~$\Mm_S^o$ be the set of~$m\in\Mm_S$ for
which~$S$ contains all primes~$p$ with~$m_p=m$, and
put~$S^o=S_{\Mm_S^o}$; this is the largest subset of~$S$ of the
form~$S_\Mm$. Similarly, put~$\widebar S=S_{\Mm_S}$, the smallest
superset of~$S$ of the form~$S_\Mm$. The techniques here can be
applied to the sets~$S^o$ and~$\widebar S$ and, since
\[
\mertens_{\widebar S}(N) \le \mertens_S(N)\le\mertens_{S^o}(N),
\]
we get some information on the asymptotic behaviour of~$\mertens_S(N)$.
\end{remarks}

The formula~\eqref{eqn:howhorrible} is particularly simple in the case
that~$\Mm$ is closed under multiplication by~$\NN$: that is, if~$a\in
\Mm$ and~$b\in\NN$, then~$ab\in\Mm$.
In this case~$\Mm$ is closed under least common multiples and,
for~$n\notin\Mm$, we have~$\bar m_n=1$. Thus the inner sum
in~\eqref{eqn:howhorrible} is empty for~$\bar m\ne 1$, and
we get\mc{newforMcase4}
\begin{equation}\label{Dancing easier}
\mertens_{S_\Mm}(N)=\sum_{\genfrac{}{}{0pt}{}{n\le N}{n\notin\Mm}}\frac{1}{n}
+C_{\ref{newforMcase4}}+\bigo\left(2^{-\sqrt{N}}\right).
\end{equation}
Provided~$\mertens_{S_\Mm}(N)\to\infty$ as~$N\to\infty$, this implies
that
\[
\mertens_{S_\Mm}(N)\sim\sum_{\genfrac{}{}{0pt}{}{n\le N}{n\notin\Mm}}\frac{1}{n}.
\]
Many, though not all, of our examples will be of this form. The task
is then to choose sets~$\Mm$ which are closed under multiplication
by~$\NN$, and for which we can control the asymptotics of the sum
in~\eqref{Dancing easier}. One technique we will often use for this
is \emph{partial} (or \emph{Abel}) \emph{summation}: if we
write~$\pi_\Mm(x)=\vert\{n\le x \mid n\not\in\Mm\}\vert$ and $f$ is a
positive differentiable function on the positive reals, then
\[
\sum_{\genfrac{}{}{0pt}{}{n\le x}{n\notin\Mm}} f(n)
= \pi_\Mm(x)f(x) + \int_1^x \pi_\Mm(t)f'(t)\dee t,
\]
with the dominant term generally coming from the integral.
In several cases the asymptotics of~$\pi_\Mm(x)$ are already well
understood.

\section{Finite sets of primes}

In order to prove Theorem~\ref{thm:dense}, we need to choose
finite sets of primes~$S$ for which we can make good estimates
for the coefficient of the leading term in Mertens' Theorem.
These calculations are simplified by considering only
primes~$p$ for which~$m_p$ is prime.

Let~$\Ll$ be a finite set of primes and take~$\Mm=\Ll$, so
that
\[
S=S_\Ll=\{p\in\PP\mid m_p\in\Ll\},
\]
which is a finite set. By~\cite[Theorem~1.4]{EMSW}, we have
\[
\mertens_{S_\Ll}(N)= k_\Ll\log(N) + C_\Ll +\bigo\left(N^{-1}\right),
\]
for some~$k_\Ll\in(0,1]\cap\mathbb Q$ and constant~$C_\Ll$. The
following lemma gives upper and lower bounds for~$k_\Ll$.

\begin{lemma}\label{lem:klapprox} Let~$\Ll$ be a finite subset of~$\PP$.
\begin{enumerate}
\item We have $\displaystyle{k_\Ll\ \le\
\prod_{\ell\in\Ll}\(1-\frac 1\ell+\frac 1{\ell(2^\ell-1)}\)}$.
\item For~$\ell\in\PP\setminus\Ll$, we have~$\(1-\frac 1\ell\)k_\Ll\
\le\ k_{\Ll\cup\{\ell\}}$.
\end{enumerate}
\end{lemma}

\begin{proof} For~$\Ll'$ a subset of~$\Ll$, we
write~$m(\Ll')=\prod_{\ell\in\Ll'}\ell$. We break up the Mertens sum
as in~\eqref{eqn:howhorrible}, noting
that~$\widebar\Mm=\{m(\Ll')\mid \Ll'\subseteq\Ll\}$:
\begin{equation}\label{eqn:mertensfinite}
\mertens_{S_\Ll}(N)\sim
\sum_{\Ll'\subseteq\Ll}\frac{\lvert2^{m(\Ll')}-1\rvert_{S_{\Ll'}}}{m(\Ll')}
\sum_{\genfrac{}{}{0pt}{}{n\le N/m(\Ll')}
{\ell\notdivides n\text{ for }\ell\in\Ll\setminus\Ll'}}
\frac{\lvert n\rvert_{S_{\Ll'}}}{n}.
\end{equation}
By~\cite[Proposition~5.2]{EMSW}, we have
\[
\sum_{n\in N}\frac{\lvert n\rvert_{S_{\Ll'}}}{n}
= k'_{\Ll'}\log N+C'_{\Ll'}+\bigo\left(N^{-1}\right),
\]
with $k'_{\Ll'}=\prod_{p\in S_{\Ll'}}\frac p{p+1}$. Moreover,
by~\cite[Lemma~5.1]{EMSW},
\[
\sum_{\genfrac{}{}{0pt}{}{n\le N}{\ell\notdivides n\text{ for }\ell\in\Ll\setminus\Ll'}}
\frac{\lvert n\rvert_{S_{\Ll'}}}{n}
= k'_{\Ll'}\prod_{\ell\in\Ll\setminus\Ll'}
\left(1-\frac{\lvert \ell\rvert_{S_{\Ll'}}}{\ell}\right) \log N
+ C''_{\Ll'}+\bigo\left(N^{-1}\right).
\]
In particular, the coefficient of the~$\log N$ term
is~$\prod_{\ell\in\Ll\setminus\Ll'}\left(1-\frac
1\ell\right)\prod_{p\in S_{\Ll'}\setminus\Ll}\frac p{p+1}$, which is
at most~$\prod_{\ell\in\Ll\setminus\Ll'}\left(1-\frac
1\ell\right)$. Moreover, for~$\ell\in\Ll'$ and~$p$ such
that~$m_p=\ell$, we have $\ord_p(2^{m(\Ll')}-1)\ge\ord_p(2^\ell-1)$ so
\[
\lvert2^{m(\Ll')}-1\rvert_{S_{\Ll'}}\le \prod_{\ell\in\Ll'}\frac 1{2^\ell-1}.
\]
Putting everything back into~\eqref{eqn:mertensfinite}, we see that
the coefficient~$k_\Ll$ of the~$\log N$ term is bounded above by
\[
\sum_{\Ll'\subseteq\Ll} \prod_{\ell\in\Ll'} \frac 1{\ell(2^{\ell-1})}
\prod_{\ell\in\Ll\setminus\Ll'}\left(1-\frac 1\ell\right)
= \prod_{\ell\in\Ll} \left(1-\frac 1\ell+\frac 1{\ell(2^{\ell-1})}\right).
\]
This proves (i), and the proof of (ii) is similar but easier: we have
\[
\mertens_{S_{\Ll\cup\{\ell\}}}(N)\sim \sum_{\Ll'\subseteq\Ll\cup\{\ell\}}
\frac{\lvert2^{m(\Ll')}-1\rvert_{S_{\Ll'}}}{m(\Ll')}
\sum_{\genfrac{}{}{0pt}{}{n\le N/m(\Ll')}{\gcd(n,\ell m(\Ll))\divides m(\Ll')}}
\frac{\lvert n\rvert_{S_{\Ll'}}}{n},
\]
and, for~$\Ll'$ a subset contained in~$\Ll$, the contribution of the
sum corresponding to~$\Ll'$ is~$\left(1-\frac{\lvert
  \ell\rvert_{S_{\Ll'}}}{\ell}\right)$ times the contribution of the
sum corresponding to~$\Ll'$ in~\eqref{eqn:mertensfinite}. In
particular, the Mertens sum for~$\Ll\cup\{\ell\}$ is at
least~$\left(1-\frac 1\ell\right)$ times that for~$\Ll$.
\end{proof}

\begin{proof}[Proof of Theorem~\ref{thm:dense}]
Let~$k\in(0,1)$ and~$\e>0$, and choose two
primes~$\ell_0>1+\frac{k}\e$ and~$\ell_1>\ell_0$ such
that~$\prod_{\ell_0\le\ell<\ell_1}\(1-\frac 1l+\frac 1{l(2^l-1)}\)<
k$; this is possible since the product over all primes greater
than~$\ell_0$ converges to~$0$.

We choose recursively a subset~$\Ll$
of~$\{\ell\in\PP\mid\ell<\ell_1\}$, using the greedy algorithm as
follows. Let~$\ell\in\PP$ and suppose we have already
defined~$\Ll(\ell):=\Ll\cap\{1,\ldots,\ell-1\}$. If~$k\le
k_{\Ll(\ell)}<k+\e$ then we are done and~$\Ll=\Ll(\ell)$;
otherwise~$\ell\in\Ll$ if and only if~$k_{\Ll(\ell)\cup\{\ell\}}\ge k$.

The claim is then that, for the subset~$\Ll$ given by this algorithm,
the leading coefficient~$k_\Ll$ satisfies~$k\le k_{\Ll}<k+\e$. The
first inequality is clear from the definition, while the second
follows from the following two observations:
\begin{enumerate}
\item There is a prime~$\ell$ with $\ell_0\le\ell<\ell_1$ such
that~$\ell\notin\Ll$: if not, by Lemma~\ref{lem:klapprox}(i),
\[
k_\Ll\le\prod_{\ell_0\le\ell<\ell_1}\(1-\frac 1\ell+\frac 1{\ell(2^\ell-1)}\)<k,
\]
which is absurd.
\item With $\ell$ as in (i), we have~$k_{\Ll(\ell)\cup\{\ell\}}<k$,
since~$\ell\notin\Ll$; thus, by Lemma~\ref{lem:klapprox}(ii),
\[
k_{\Ll(\ell)}\le \left(\frac \ell{\ell-1}\right)k_{\Ll(\ell)\cup\{\ell\}}
< \left(\frac {\ell_0}{\ell_0-1}\right)k<k+\e.
\]
\end{enumerate}
\end{proof}

\begin{remark}\label{rmk:sublog}
Let~$\Ll$ be an \emph{infinite} set of primes such
that~$\sum_{\ell\in\Ll}\frac 1\ell$ diverges and
put~$S=S_{\Ll}=\{p\in\PP\mid m_p\in\Ll\}$.
Then~$\mertens_S(N)\le\mertens_{S_{\Ll'}}(N)$, for any finite
subset~$\Ll'$ of~$\Ll$, so~$\mertens_S(N)$ grows no more
quickly than~$k_{\Ll'}\log N$. Since~$k_{\Ll'}$ is at
most~$\prod_{\ell\in\Ll'}\(1-\frac 1\ell+\frac
1{\ell(2^\ell-1)}\)$, by Lemma~\ref{lem:klapprox}(i), there are
finite subsets~$\Ll'$ of~$\Ll$ with~$k_{\Ll'}$ arbitrarily
close to~$0$, and we deduce
that~$\mertens_S(N)=\littleo\left(\log N\right)$.
\end{remark}

\section{Logarithmic growth for infinite sets of primes}
\label{S:logarithmic}

In this section we will prove Theorem~\ref{thm:onto}.
Fix~$\ell\in\PP$, and let~$\widebar\Mm_{\ell}=
\{n\in\mathbb N\mid\ell\divides n\}$, so that
\[
S_{\widebar\Mm_{\ell}}=\widebar S_{\ell}=\{p\in\PP\mid \ell\divides m_p\}.
\]
By Hasse~\cite{MR0186653,MR0205975} these sets have a positive
Dirichlet density within the set of primes: for~$\ell>2$ the density
is~$\frac{\ell}{\ell^2-1}$, and for~$\ell=2$ the density
is~$\frac{17}{24}$. They also have natural density by, for
example,~\cite[Theorem~2]{MR736730}. Noting that~$\widebar\Mm_\ell$
is closed under multiplication by~$\NN$, by~\eqref{Dancing easier} we
have\mc{ellconstantA}\mc{ellconstantB}
\begin{align*}
\mertens_{\widebar S_{\ell}}(N)&=
\sum_{\genfrac{}{}{0pt}{}{n\le N}{\ell\notdivides n}}\frac1n
+ C_{\ref{ellconstantB}}+\bigo\left(2^{-\sqrt{N}}\right) \\
&=\(1-\textstyle\frac{1}{\ell}\)\log N
+C_{\ref{ellconstantA}}+\bigo\left(N^{-1}\right).
\end{align*}
Indeed, if~$\Ll$ is any \emph{finite} set of primes, applying the same
argument to
\[
\widebar \Mm_{\mathcal{L}}=\{n\in\mathbb N\mid \ell\divides n
\mbox{ for some }\ell\in\mathcal{L}\},
\]
and~$\widebar S_{\mathcal{L}}=
\{p\in\PP\mid m_p\in \widebar M_{\mathcal{L}}\}$, we get\mc{ellconstantA'}
\[
\mertens_{\widebar S_{\Ll}}(N)=
\prod_{\ell\in\Ll}\(1-\textstyle\frac{1}{\ell}\)\log N
+C_{\ref{ellconstantA'}}+\bigo\left(N^{-1}\right).
\]
Since the set
$\left\{\prod_{\ell\in\Ll}\(1-\textstyle\frac{1}{\ell}\)\mid
\Ll\subset\PP\mbox{ finite}\right\}$ is dense in~$[0,1]$, this
gives an easy way of
getting a dense set of values for the leading coefficient in Mertens'
Theorem. Note however that Theorem~\ref{thm:dense} was more delicate,
since the claim was that a dense set of values can be obtained using
only~\emph{finite} sets $S$. Similarly, Theorem~\ref{thm:onto} claims
more: \emph{every} value in~$(0,1)$ can be obtained as leading
coefficient.

\begin{proof}[Proof of Theorem~\ref{thm:onto}.]
Now let~$\mathcal{L}\subset\PP$ be any set of primes for which the
product~$k_{\mathcal{L}}:=
\prod_{\ell\in\Ll}\(1-\textstyle\frac{1}{\ell}\)$
is non-zero, define~$\widebar \Mm_\Ll$ and~$\widebar S_\Ll$ as above,
and apply the argument above to obtain
\[
\mertens_{\widebar S_{\mathcal{L}}}(N)\sim
\sum_{\genfrac{}{}{0pt}{}{n\le N}{n\notin\widebar \Mm_{\mathcal{L}}}}
\frac 1n.
\]
Applying~\cite[Theorem~I.3.11]{MR1342300} we have that
$\vert\{n\le x\mid n\notin\widebar \Mm_\Ll\}\vert\sim k_\Ll x$ thus, by
partial summation, we get
\[
\mertens_{\widebar S_{\mathcal{L}}}(N)
\sim
k_{\mathcal{L}}\log N.
\]
This gives Theorem~\ref{thm:onto} since
$\left\{k_\Ll\mid\mathcal{L}\subseteq\PP\right\}=[0,1]$.
\end{proof}

\section{Sublogarithmic growth}

Now we consider sets giving intermediate sublogarithmic
growth, proving Theorem~\ref{thm:logdelta}. We return to sets close
to~$\widebar \Mm_\Ll$ and~$\widebar S_\Ll$ of~\S\ref{S:logarithmic}
but now for infinite sets of primes~$\Ll$ such
that~$\prod_{\ell\in\Ll}\(1-\frac 1\ell\)=0$. We will need a result from
analytic number theory that allow sets of primes to be selected with
prescribed properties, whose proof we defer to~\S\ref{S:proofofL}.

\begin{proposition}\label{justaperfectday}
For any~$\d\in(0,1]$, there is a set of primes~$\Ll$ such that
\begin{equation}\label{eqn:A15}
\sum_{\genfrac{}{}{0pt}{}{\ell\le x}{l\in\Ll}}\frac{\log\ell}{\ell} =
\delta\log x+\bigo(1)
\end{equation}
and, for any~$c>1$, there is a set of primes~$\Ll'\subseteq\Ll$ such that
\[
\prod_{p\in\Ll'}\left(1+\frac{1}{p}\right)=c
\qquad\hbox{ and }\qquad
\sum_{p\in\Ll'}\frac{\log p}{p}<\infty.
\]
\end{proposition}

\begin{proof}[Proof of Theorem~\ref{thm:logdelta}]
Let~$\delta\in(0,1]$ and $k>0$, and let $\Ll$ be a set of primes
satisfying~\eqref{eqn:A15}. As before we
put~$\widebar\Mm_\Ll=\{n\in\NN\mid\ell\divides n\hbox{ for some
}\ell\in\Ll\}$ and now we set
\[
\widebar\Mm'_\Ll = \{n\in\NN\mid n\in\widebar\Mm_\Ll
\hbox{ or $n$ is not square-free}\}
\]
and $\widebar S'_\Ll=\{p\in\PP\mid m_p\in\widebar\Mm'_\Ll\}$.
Note that~$\widebar\Mm'_\Ll$ is also closed under
multiplication by~$\NN$ so that\mc{ellconstantsub}
\[
\mertens_{\widebar S'_{\Ll}}(N) =
\sum_{\genfrac{}{}{0pt}{}{n\le N}{n\not\in\widebar\Mm'_\Ll}}\frac1n
+ C_{\ref{ellconstantsub}}+\bigo\left(2^{-\sqrt{N}}\right)
\]
by~\eqref{Dancing easier}. Now we
apply~\cite[Theorem~A.5]{MR2647984} with, in the notation used
there, the function
\[
g(n) = \begin{cases} \frac 1n &\hbox{ if }n\notin\Mm'_\Ll, \\
0&\hbox{ otherwise.} \end{cases}
\]
Note that, by~\eqref{eqn:A15}, the
hypotheses~\cite[(A.15--17)]{MR2647984} of that Theorem are indeed
satisfied. We conclude that
\[
\sum_{\genfrac{}{}{0pt}{}{n\le N}{n\not\in\widebar\Mm'_\Ll}}\frac1n =
k_\Ll(\log N)^\delta + \bigo\((\log N)^{\delta-1}\),
\]
where~$k_\Ll>0$ is
\[
k_\Ll = \frac 1{\Gamma(\delta+1)}
\prod_{p\in\PP}\(1-\frac 1p\)^\delta\prod_{p\notin\Ll}\(1+\frac 1p\)
\]
by~\cite[(A.24)]{MR2647984}. Notice that we can adjust~$\Ll$ by
any set of primes~$\Ll'$ such
that~$\sum_{\ell\in\Ll'}\frac{\log\ell}{\ell}$ converges
without affecting the hypothesis~\eqref{eqn:A15}.

Assume now that~$\Ll$ is the set of primes~$\Ll$ constructed in
Proposition~\ref{justaperfectday}.
Let~$\Ll''\subseteq\PP\setminus\Ll$ be a set of primes such
that~$k\prod_{p\in\Ll''}\(1+\frac 1p\)\ge k_\Ll$. By
Proposition~\ref{justaperfectday}, there is a subset~$\Ll'$
of~$\Ll$ such that~$\prod_{p\in\Ll'}\(1+\frac 1p\)=
k/k_\Ll\prod_{p\in\Ll''}\(1+\frac 1p\)$. In particular,
putting~$\Ll_0=(\Ll\setminus\Ll')\cup\Ll''$, we
have~$k_{\Ll_0}=k$ so the set~$\widebar S'_{\Ll_0}$ gives the
required asymptotic.
\end{proof}

\begin{remark}
Since the sets~$\widebar\Mm_\Ll$ and~$\widebar\Mm'_\Ll$ coincide on the
set of square-free natural numbers, there is a constant $c_\Ll$ such
that
\[
\sum_{\genfrac{}{}{0pt}{}{n\le N}{n\notin\widebar\Mm'_{\mathcal{L}}}}\frac1n \le
\sum_{\genfrac{}{}{0pt}{}{n\le N}{n\notin\widebar\Mm_{\mathcal{L}}}}\frac1n \le
c_\Ll \sum_{\genfrac{}{}{0pt}{}{n\le N}{n\notin\widebar\Mm_{\mathcal{L}}}}\frac1n.
\]
In particular, the Mertens sum~$\mertens_{\widebar S_\Ll}(N)$ also grows
like~$\(\log N\)^\delta$. We have chosen to use~$\widebar\Mm'_\Ll$ here
rather than~$\widebar\Mm_\Ll$ since it is for such a set that we were able
to find an off-the-shelf reference~\cite[Theorem~A.5]{MR2647984} for the
asymptotics.
\end{remark}

\section{Doubly logarithmic growth}

Here we consider sets giving doubly logarithmic growth or slower, in
particular proving Theorem~\ref{thm:loglog}. In the case~$r=1$, the
proof is based on the taking the set~$\Mm$
of~\S\ref{In a new world, me and you, girl} to be the
set~$\NN\setminus\PP$ of composite natural numbers, so that~$S=S_\Mm$
is the set of primes~$p$ such that~$m_p$ is composite. Since~$\Mm$ is
closed under multiplication by~$\NN$, applying~\eqref{Dancing easier}
we have\mc{mertensprimes}
\[
\mertens_{S}(N)=\sum_{\genfrac{}{}{0pt}{}{p\le N}{p\in\PP}}\frac{1}{p}
+C_{\ref{mertensprimes}}+\bigo\left(2^{-\sqrt{N}}\right).
\]
By Mertens' original theorem~\cite{mertens}, we
have\mc{mertensclassical}
\begin{equation}\label{eqn:mertens}
\sum_{\genfrac{}{}{0pt}{}{p\le N}{p\in\PP}}\frac{1}{p}
=\log\log N+C_{\ref{mertensclassical}}+\bigo\left((\log N)^{-1}\right),
\end{equation}
and hence\mc{mertensprimes2}
\[
\mertens_{S}(N)=\log\log N+C_{\ref{mertensprimes2}}
+\bigo\left((\log N)^{-1}\right),
\]
which is an improved form (i.e. with error term) of
Theorem~\ref{thm:loglog} with~$k=1$ and~$r=1$.

\begin{remark}
The complement of this set~$S$ is the set of primes~$p$ for
which~$m_p$ is prime
so~$\mertens_{\PP\setminus S}(N)=\littleo\left(\log N\right)$, by
Remark~\ref{rmk:sublog}. Thus both $\mertens_{S}(N)$ and
$\mertens_{\PP\setminus S}(N)$ are~$\littleo\left(\log N\right)$.
\end{remark}

For the general case of Theorem~\ref{thm:loglog} we will need the
following lemma, which gives asymptotics for the number of integers
with exactly~$r$ prime factors (counted with multiplicity), all from
a fixed set of primes. For~$\Ll$ a set of primes and~$n\in\NN$ we denote
by~$\Omega_\Ll(n)$ the number of primes factors of~$n$ in~$\Ll$
(counted with multiplicity), and abbreviate~$\Omega(n)=\Omega_\PP(n)$.

\begin{lemma}
Let~$\Ll$ be a set of primes of natural density~$\d$ and $r\in\NN$. Then
\begin{equation}\label{eqn:piLr}
\vert\left\{n\le x\mid \Omega_\Ll(n)=\Omega(n)=r\right\}\vert
\sim \d^r \frac x{\log x}\frac{\(\log\log x\)^{r-1}}{(r-1)!}.
\end{equation}
\end{lemma}

\begin{proof}
When~$r=1$, the case~$\Ll=\PP$ is the prime number theorem and the case
of general~$\Ll$ follows immediately, since~$\Ll$ has density~$\d$.
For~$r>1$, the case~$\Ll=\PP$ is a result of
Landau~\cite[XIII~\S56~(5)]{MR0068565d}, proved by induction on~$r$.
The proof (using the prime number theorem and partial summation --
see~\cite[\S7.4]{MR2378655} for a sketch) works
equally well for any set~$\Ll$ as in the lemma, and the result follows.
\end{proof}

\begin{proof}[Proof of Theorem~\ref{thm:loglog}]
Let~$r\in\NN$ and~$k>0$. We pick a natural number~$m$ such
that
\[
k_m:=\sum_{d\divides m}\frac 1d > k(r!)
\]
and a set~$\Ll$ of primes with natural
density~$\d=\(k(r!)/k_m\)^{1/r}$. Denote by~$\Mm_{r,\Ll,m}$ the
set of natural numbers~$n$ such that
either~$\Omega(n/\gcd(m,n))>r$ or~$n/\gcd(m,n)$ has a prime
factor outside~$\Ll$.

We put $S=S_{\Mm_{r,\Ll,m}}$ and apply~\eqref{Dancing easier} to get
\begin{equation}\label{eqn:breakup}
\mertens_{S}(N)\sim
\sum_{\genfrac{}{}{0pt}{}{n\le N}{n\notin\Mm_{r,\Ll,m}}}\frac{1}{n}
= \sum_{d\divides m} \frac 1d
\sum_{\genfrac{}{}{0pt}{}{n\le N/d}{n\notin\Mm_{r,\Ll,1}}}\frac{1}{n}.
\end{equation}
On the other hand, by~\eqref{eqn:piLr}, we have
\[
\pi_{\Mm_{r,\Ll,1}}(x) =
\vert\left\{n\le x\mid \Omega(n)=\Omega_\Ll(n)\le r\right\}\vert
\sim \d^r \frac x{\log x}\frac{\(\log\log x\)^{r-1}}{(r-1)!}.
\]
Applying partial summation gives
\[
\sum_{\genfrac{}{}{0pt}{}{n\le x}{n\notin\Mm_{r,\Ll,1}}}\frac{1}{n}
\sim  \d^r \frac{\(\log\log x\)^r}{r!},
\]
and substituting this into~\eqref{eqn:breakup} gives the result,
because of the choice of~$\d$.
\end{proof}

\begin{remark}
Let $\theta$ be any positive, increasing,
differentiable function on the positive reals such that, for large
enough $x$, both $\theta(x)\le\log\log x$ and $\theta'(x)\le\frac
1{x\log x}$. Then there is a set of primes $\Ll_\theta$ such that
\[
\sum_{\genfrac{}{}{0pt}{}{p\le x}{p\in\Ll_\th}}\frac 1p\sim\theta(x)
\]
and hence, putting $S_\theta=\{p\in\PP\mid m_p\not\in\Ll_\theta\}$ and
applying~\eqref{Dancing easier}, we have
\[
\mertens_{S_\theta}(N)\sim\theta(x).
\]
The existence of such a set of primes~$\Ll_\theta$ comes
from~\eqref{eqn:mertens} and the following lemma, whose proof using
the greedy algorithm is straightforward but technical so is omitted. It
seems almost certain that a lemma of this sort exists in the literature
but we have not been able to find it.
\end{remark}

\begin{lemma} Suppose~$f$ is a positive, increasing, differentiable
function on the positive reals and~$a_n$ are non-negative
reals converging to $0$ such that~$\sum_{n\le x}a_n\sim f(x)$. We
write~$\d(x)=\sum_{n\le x}a_n- f(x)$. Suppose we have a positive,
increasing, differentiable function~$\th$ on the positive reals such that:
\begin{enumerate}
\item there is an~$x_0>0$ such that~$\th(x)\le f(x)$ and ~$\th'(x)\le
  f'(x)$ for all~$x>x_0$;
\item $\d=\littleo(\th)$.
\end{enumerate}
Then there is a subset~$\Nn_\th\subset\NN$ such that
\[
\sum_{\genfrac{}{}{0pt}{}{n\le x}{n\in\Nn_\th}}a_n\ \sim\ \th(x).
\]
\end{lemma}

\section{Convergence for co-infinite sets of primes}

\begin{proof}[Proof of Proposition~\ref{prop:zero}.]
Fix a prime~$\ell$ and set~$\widebar\Mm_\ell^c=\{n\in\mathbb
N\mid\ell\notdivides n\}$, the complement of the set~$\widebar\Mm_\ell$
considered in~\S\ref{S:logarithmic};
thus~$S=S_{\widebar\Mm_\ell^c}=\{p\in\PP\mid\ell\notdivides m_p\}$
is a set of primes with positive natural density.
Although~$\widebar\Mm_\ell^c$ is not closed under multiplication
by~$\NN$, it is closed under least common multiples; moreover,
for~$m\in\widebar\Mm_\ell^c$, we have $\NN_m=\{m\ell^e\mid
e\ge 0\}$, in the notation of~\S\ref{In a new world, me and you,
  girl}. Thus, from~\eqref{eqn:howhorrible},\mc{newforzero}
\begin{equation}\label{eqn:zero}
\mertens_{S}(N)=\sum_{\genfrac{}{}{0pt}{}{m\in\NN}{\ell\notdivides m}}
\frac{\lvert2^m-1\rvert_{S_m}}m
\sum_{1< \ell^e\le N/m} \frac{\lvert\ell\rvert_{S_m}^e}{\ell^e}
+C_{\ref{newforzero}}+\bigo\left(2^{-\sqrt{N}}\right).
\end{equation}
Now\[
\sum_{1< \ell^e\le N/m}
\frac{\lvert\ell\rvert_{S_m}^e}{\ell^e} \le
\sum_{e\ge 1} \frac 1{\ell^e} =
\frac 1{(\ell-1)}.
\]
Thus the terms of the (outer) sum in~\eqref{eqn:zero} converge and,
since $\lvert 2^m-1\rvert_{S_m}=\(2^m-1\)^{-1}$, the difference
between each term and its limit is\mc{boundzero}
\[
\frac{\lvert2^m-1\rvert_{S_m}}m \sum_{\ell^e> N/m}
\frac{\lvert\ell\rvert_{S_m}^e}{\ell^e} \le
\frac m{2^m-1}\sum_{e>\frac{\log(N/m)}{\log\ell}} \frac 1{\ell^e}
\le C_{\ref{boundzero}} \frac 1{2^mN}.
\]
Plugging this back into the sum in~\eqref{eqn:zero}, we see that it
converges and the difference between it and its limit is bounded by
\[
\frac{C_{\ref{boundzero}}}N\sum_{\genfrac{}{}{0pt}{}{m\in\NN}{\ell\notdivides m}}
\frac 1{2^m} = \bigo\left(N^{-1}\right).
\]
\end{proof}

\begin{remark}
It is straightforward to generalize this proof to the case where~$\Mm$
is the complement of the set~$\widebar\Mm_\Ll$ considered
in~\S\ref{S:logarithmic}, for any finite set of primes~$\Ll$, so
that~$S_\Mm$ is the set of primes~$p$ such that~$m_p$ is not divisible
by any~$\ell\in\Ll$. By a special case of a very general result of
Wiertelak~\cite[Theorem~2]{MR736730}, when~$\Ll$ consists only of odd
primes the set~$S_\Mm$ has natural
density~$\prod_{\ell\in\Ll}\(1-\frac{\ell}{\ell^2-1}\)$. In
particular, this density can be arbitrarily close to~$0$.
\end{remark}

\section{Transcendental constants}

Our first example of a transcendental constant comes from an
elementary result in analytic number theory. Let~$\Mm$ be the set of
non-squarefree natural numbers; then a theorem of Landau gives
\[
\pi_\Mm(x) = \vert\{n\le x\mid n\notin\Mm\}\vert = \frac 6{\pi^2}x+\littleo(\sqrt x)
\]
(see for example~\cite[XLIV~\S162]{MR0068565d}
or~\cite[Theorem~I.3.10]{MR1342300}). Thus, by partial summation
and~\eqref{Dancing easier}, we get\mc{sqfree}
\[
\mertens_{S_\Mm}(N) = \frac
6{\pi^2} \log N + C_{\ref{sqfree}} + \littleo\(N^{-1/2}\).
\]

For our second example, fix a prime~$\ell$ and
set~$\Mm_{(\ell)}=\{\ell^e\mid e\ge 0\}$, so
that~$S=S_{\Mm_{(\ell)}}$ is the infinite set of primes~$p$ for
which~$m_p$ is a power of~$\ell$. This is a thin set of primes: that
is, it has density zero.
As in the previous section, the set~$\Mm_{(\ell)}$ is closed under
least common multiples, but not under multiplication
by~$\NN$. Applying~\eqref{eqn:howhorrible}, we get\mc{newfortrans}
\begin{equation}\label{eqn:trans}
\mertens_{S}(N)=\sum_{e=0}^{\infty} \frac 1{\ell^e}
\sum_{\genfrac{}{}{0pt}{}{2\le n\le N}{\ord_\ell(n)=e}}
\frac{\lvert2^{n}-1\rvert_{S_e}}{n}+C_{\ref{newfortrans}}
+\bigo\left(2^{-\sqrt{N}}\right),
\end{equation}
where~$S_e$ is the finite set of primes dividing~$2^{\ell^e}-1$. (This
set was denoted by~$S_{\ell^e}$ in~\eqref{eqn:howhorrible}.)
Noting that~$\ell\notin S$, we observe that, for any~$e\ge
0$,~$n\in\NN$ such that~$\ord_\ell(n)=e$, and prime~$p$
dividing~$2^{\ell^e}-1$, by~\eqref{eqn:ordp2n-1} we have
\[
\ord_p\left(2^{n}-1\right) =
\ord_p\left(2^{\ell^e}-1\right)+\ord_p(n).
\]
Since every prime divisor of~$2^{\ell^e}-1$ lies in~$S_e$, we deduce that
\[
\lvert 2^n-1\rvert_{S_e} = \frac {\lvert n\rvert_{S_e}}{2^{\ell^e}-1}.
\]
Thus the sum in~\eqref{eqn:trans} becomes
\begin{equation}\label{eqn:trans2}
\sum_{e=0}^{\infty} \frac 1{\ell^e(2^{\ell^e}-1)}
\sum_{\genfrac{}{}{0pt}{}{2\le n\le N/\ell^e}{\ell\notdivides n}}
\frac{\lvert n\rvert_{S_e}}{n}.
\end{equation}
Now, by~\cite[Proposition~5.2, Lemma~5.1]{EMSW}, we have
\begin{equation}\label{eqn:bigoe}
\sum_{\genfrac{}{}{0pt}{}{2\le n\le N}{\ell\notdivides n}}
\frac{\lvert n\rvert_{S_e}}{n} =
\left(1-\frac 1\ell\right) k_e\log N+ \bigo_e\left(1\right),
\end{equation}
where~$k_e=\prod_{p\in S_e}\frac p{p+1}$. Here we need to control the
error terms uniformly in~$e$. For this, we use the following lemma,
which we will prove at the end of the section.

\begin{lemma}\label{lem:bounderror}
For~$S'$ any finite set of primes
put~$k'_{S'}=\prod_{p\in S'}\frac p{p+1}$ and
\[
f_{S'}(N)=\sum_{n\le N}\frac{\lvert n\rvert_{S'}}{n} - k'_{S'}\log N.
\]
By~\cite[Proposition~5.2]{EMSW}, there exists~$A_{S'}>4$ such
that~$\lvert f_{S'}(N)\rvert\le A_{S'}$, for all~$N>1$.

Now fix~$S'$, let~$p\in\PP\setminus S'$ and put~$S''=S'\cup\{p\}$.
Then~$\lvert f_{S''}(N)\rvert\le 2A_{S'}$, for all~$N>1$.
\end{lemma}

In particular the~$\bigo_e(1)$ error in~\eqref{eqn:bigoe}
is~$\bigo\(2^{\vert S_e\vert}\)$, with an implied constant
independent of~$e$, and~$2^{\vert S_e\vert}\le \prod_{p\in S_e}
p \le 2^{\ell^e}-1$. Thus the error in each term of the outside
sum in~\eqref{eqn:trans2} is~$\bigo(1/\ell^e)$ and the sum of
these errors converges. Thus~\eqref{eqn:trans}
and~\eqref{eqn:trans2} give
\[
\mertens_S(N)\sim k_S\log N,
\]
with
\[
k_S = \sum_{e=0}^{\infty} \frac {(\ell -1)}{\ell^{e+1}(2^{\ell^e}-1)}
\prod_{p\in S_e}\frac p{p+1}.
\]
Now the partial sums give infinitely many rational
approximations~$\frac ab$ of~$k_S$ with error~$\bigo\(b^{-\ell}\)$;
thus, provided~$\ell\ge 3$, we deduce that~$k_S$ is transcendental by
Roth's Theorem.

It only remains to prove Lemma~\ref{lem:bounderror}.

\begin{proof}[Proof of Lemma~\ref{lem:bounderror}]
We have
\[
\sum_{n\le N}\frac{\lvert n\rvert_{S''}}n =
\sum_{r=0}^{\lfloor\log N/\log p\rfloor} \frac 1{p^{2r}}
\sum_{\genfrac{}{}{0pt}{}{n\le N/p^r}{p\notdivides n}}
\frac{\lvert n\rvert_{S'}}n
\]
and
\[
\sum_{\genfrac{}{}{0pt}{}{n\le N}{p\notdivides n}}
\frac{\lvert n\rvert_{S'}}n =
\(1-\frac 1p\)k'_{S'}\log N - \frac{k'_{S'}\log p}{p}
+ f_{S'}(N) -\frac 1p f_{S'}(N/p).
\]
Putting these together, we get
\begin{eqnarray*}
f_{S''}(N)&=&\frac{k'_{S'}(p-1)}p
\sum_{r>\left\lfloor\frac{\log N}{\log p}\right\rfloor}\frac 1{p^{2r}} \log N
-\frac{k'_{S'}(p-1)\log p}p
\sum_{r=0}^{\left\lfloor\frac{\log N}{\log p}\right\rfloor} \frac r{p^{2r}}\\[5pt]
&&-\frac{k'_{S'}\log p}{p}
\sum_{r=0}^{\left\lfloor\frac{\log N}{\log p}\right\rfloor} \frac 1{p^{2r}}
+\sum_{r=0}^{\left\lfloor\frac{\log N}{\log p}\right\rfloor}\frac 1{p^{2r}}
\(f_{S'}\(N/p^r\)-\frac 1p f_{S'}\(N/p^{r+1}\)\).
\end{eqnarray*}
Using~$0<k'_{S'}\le 1$,~$p\ge 3$ and~$N\ge 2$, the first three
terms are absolutely bounded by~$\frac{pk'_{S'}}{(p+1)\log N}
<\frac 1{\log 2}$, $\frac{pk'_{S'}}{(p+1)(p^2-1)}<\frac 3{32}$
and $\frac{pk'_{S'}\log p}{(p^2-1)}<\frac {3\log 3}{8}$
respectively, whose sum is bounded by~$2$. The final term is
bounded in absolute value by~$\frac p{p-1}A_{S'}<\frac 32
A_{S'}$ and the result follows from the assumption
that~$A_{S'}>4$.
\end{proof}

\begin{remark}
In fact~\cite[Proposition~5.2]{EMSW} says
that~$f_{S'}(N)=C'_{S'}+\bigo\(N^{-1}\)$ so a finer analysis of the errors
in~\eqref{eqn:bigoe} along the lines of Lemma~\ref{lem:bounderror} should
allow one to get an asymptotic expression for~$M_S(N)$ with an error term.
\end{remark}

\section{Existence of suitable sets of primes}
\label{S:proofofL}

It remains only to prove Proposition~\ref{justaperfectday}.

\begin{proof}
Let~$\d\in(0,1]$. We seek first a set of primes~$\Ll$ such that
\[
\sum_{\genfrac{}{}{0pt}{}{\ell\le x}{l\in\Ll}}\frac{\log\ell}{\ell} =
\delta\log x+\bigo(1).
\]
For rational~$\delta$, such a set exists from Dirichlet's
Theorem on primes in arithmetic progression
(see~\cite[Theorem~7.3]{MR0434929}); thus~$\Ll$ would be a set of
primes defined by congruence conditions and~$\delta$ would in fact be
the natural density of~$\Ll$.
For arbitrary~$\delta$ a more delicate construction is needed.
Let~$S$ be the set of primes
in the union of intervals
\[
\bigcup_{n\in\NN}\(2^{n},2^{n+\d}\right].
\]
Now the prime number theorem implies that
\[
\pi_{\PP}(x):=\vert\{p\le x\mid p\in\PP\}\vert = \frac{x}{\log x}
+\frac{x}{(\log x)^2} +2\frac{x}{(\log x)^3}
+\bigo\left(\frac{x}{(\log x)^4}\right)
\]
and applying partial summation gives
\[
\sum_{2^{n}<p\le 2^{n+\d}}\frac{\log p}p =
\d \log 2+\bigo\(\frac 1{n^2}\).
\]
Summing over~$n$ gives the required asymptotic.

For~$\Ll'\subseteq\Ll$,
write~$\Sigma(\Ll')=\sum_{p\in\Ll'}\textstyle\frac{\log p}{p}$.
Taking logarithms, the statement now sought is that any~$a>0$
can be written as
\[
\left.
\begin{array}{l}
a=\displaystyle\sum_{p\in\Ll'}\log\(
1+\textstyle\frac{1}{p}
\),\\[15pt]
\Sigma(\Ll')<\infty.
\end{array}\right\}\label{dranksangriainthepark}
\]
The basic idea is to use the greedy algorithm but on a subset of~$\Ll$
which is forced to be sparse enough to ensure the convergence
of~$\Sigma(\Ll')$.

Let
\begin{equation} \label{L'}
\Ll':=\mathbb{P}\cap
\left(
\bigcup_{X<m< Y} \left(\left. 2^m,2^{m+\delta}\right]\right.
\cup \left(\left.2^{Y+\delta/4},2^{Y+\delta}\right]\right. \cup
\bigcup_{n=1}^{\infty}\left(\left. 2^{R_n}, 2^{R_nr_n}\right] \right.
\right),
\end{equation}
where
\[
R_n:=\left[2^{n/2}Y\right].
\]
Here we assume that~$X,Y\in \NN$ and~$r_n\in \RR$ are parameters satisfying
\[
3\le X<Y
\]
and
\[
1<r_n< 1+\frac{\delta}{R_n},
\]
for all~$n\ge 1$, so that the union on the right-hand side
of~\eqref{L'} is disjoint and~$\Ll'\subseteq \Ll$.
The construction involves choosing the parameters~$X,Y$ and~$r_n$
appropriately. More precisely, we show that, provided~$X$ is large enough,
there are choices of~$Y$ and~$r_n$ such that the corresponding set~$\Ll'$
has the required properties.

We first derive asymptotic estimates for the sums
\[
\sum_{2^m<p\le 2^{m+\delta}} \log\left(1+\frac{1}{p}\right), \quad
\sum_{2^Y<p\le 2^{Y+\delta/4}}  \log\left(1+\frac{1}{p}\right), \quad
\sum_{2^{R_n}<p\le 2^{R_nr_n}} \log\left(1+\frac{1}{p}\right).
\]
We have
\begin{equation}\label{homer}
\sum_{2^m<p\le 2^{m+\delta}} \log\left(1+\frac{1}{p}\right)
= \sum_{2^m<p\le 2^{m+\delta}} \left(\frac{1}{p}
+\bigo\left(\frac{1}{p^2}\right)\right)
= \sum_{2^m<p\le 2^{m+\delta}} \frac{1}{p}+\bigo\left(2^{-m}\right).
\end{equation}
Using the prime number theorem with error term in the well-known form
(see for example~\cite[\S4.1~Theorem~1]{MR1342300})\mc{pntc}
\[
\pi(x)=\int_{2}^{x} \frac{\dee t}{\log t}
+\bigo\left(x\exp\left(-C_{\ref{pntc}}\sqrt{\log x}\right)\right),
\]
$C_{\ref{pntc}}$ being a suitable positive constant,
and partial summation and integration, we deduce that\mc{lastc}
\begin{eqnarray}
\sum_{2^m<p\le 2^{m+\delta}} \frac{1}{p} &=&
\frac{1}{2^{m+\delta}}\cdot \(\sum_{2^m<p\le 2^{m+\delta}} 1\)
+ \int_{2^m}^{2^{m+\delta}} \frac{1}{t^2} \cdot
\(\sum_{2^m<p\le t} 1\)\dee t  \notag\\
&=& \frac{1}{2^{m+\delta}}\cdot
\int_{2^m}^{2^{m+\delta}} \frac{\dee t}{\log t}
+ \int_{2^m}^{2^{m+\delta}} \frac{1}{t^2} \cdot
\int_{2^{m}}^{t} \frac{\dee y}{\log y} \dee t
+ \bigo\(\exp\(-C_{\ref{lastc}}m^{1/2}\)\)  \notag\\
&=& \int_{2^m}^{2^{m+\delta}} \frac{\dee t}{t \log t}
+ \bigo\(\exp\(-C_{\ref{lastc}}m^{1/2}\)\) \notag\\
&=& \log\log 2^{m+\delta}- \log\log 2^m
+ \bigo\(\exp\(-C_{\ref{lastc}}m^{1/2}\)\)  \notag\\
&=& \log \(1+\frac{\delta}{m}\)
+ \bigo\(\exp\(-C_{\ref{lastc}}m^{1/2}\)\),
\label{marge}
\end{eqnarray}
for a suitable positive constant~$C_{\ref{lastc}}$. Combining~\eqref{homer}
and~\eqref{marge}, we obtain
\begin{equation}\label{bart}
\sum_{2^m<p\le 2^{m+\delta}} \log\(1+\frac{1}{p}\)
=\log \(1+\frac{\delta}{m}\) + \bigo\(\exp\(-C_{\ref{lastc}}m^{1/2}\)\).
\end{equation}
Similarly, we derive
\begin{equation} \label{bart2}
\sum_{2^{Y}<p\le 2^{Y+\delta/4}} \log\(1+\frac{1}{p}\)
=\log \(1+\frac{1}{4}\cdot \frac{\delta}{Y}\) + \bigo\(\exp\(-C_{\ref{lastc}}Y^{1/2}\)\)
\end{equation}
and
\begin{eqnarray}
b_n\ := \sum_{2^{R_n}<p\le 2^{R_nr_n}} \log\(1+\frac{1}{p}\)
&=& \log r_n +\bigo\(\exp\(-C_{\ref{lastc}}R_n^{1/2}\) \) \notag\\
&=& \log r_n +\bigo\(\exp\(-C_{\ref{lastc}}2^{n/4}Y^{1/2}\)\).
\label{lisa}
\end{eqnarray}

Now assume that~$X$ is large enough so that
\[
\sum_{2^m<p\le 2^{m+\delta}} \log\(1+\frac{1}{p}\) < a,
\]
for all~$m>X$. Let~$Y$ be the unique natural number satisfying
\[
\sum_{X<m\le Y} \sum_{2^m<p\le 2^{m+\delta}} \log\(1+\frac{1}{p}\)
< a \le
\sum_{X<m\le Y+1} \sum_{2^m<p\le 2^{m+\delta}} \log\(1+\frac{1}{p}\).
\]
Set
\begin{equation*} \label{a'def}
a':=a-\left(\sum_{X<m< Y} \sum_{2^m<p\le 2^{m+\delta}} \log\(1+\frac{1}{p}\) +
\sum_{2^{Y+\delta/4}<p\le 2^{Y+\delta}} \log\(1+\frac{1}{p}\)\right).
\end{equation*}
Using~\eqref{bart} and \eqref{bart2}, we have
\begin{eqnarray}
\frac{1}{5} \cdot \frac{\delta}{Y} &<& \log\(1+\frac{1}{4}\cdot \frac{\delta}{Y}\)+ \bigo\(\exp\(-C_{\ref{lastc}}Y^{1/2}\)\) = \sum\limits_{2^{Y} <p\le 2^{Y+\delta/4}} \log\(1+\frac{1}{p}\) \nonumber\\
& <  & a'
\le\sum\limits_{2^{Y} <p\le 2^{Y+\delta/4}} \log\(1+\frac{1}{p}\)+ \sum_{2^{Y+1}<p\le 2^{Y+1+\delta}} \log\(1+\frac{1}{p}\)  \notag\\
&=& \log\(1+\frac{1}{4}\cdot \frac{\delta}{Y}\)+ \log\(1+\frac{\delta}{Y+1}\)
+ \bigo\(\exp\(-C_{\ref{lastc}}Y^{1/2}\)\)
\notag\\
&<& \frac{5}{4}\cdot \frac{\delta}{Y}+ \bigo\(\exp\(-C_{\ref{lastc}}Y^{1/2}\)\)
 < \frac{4}{3}\cdot \log\(1+\frac{\delta}{Y}\),
\label{uter}
\end{eqnarray}
provided~$Y$ is sufficiently large (which is the case if~$X$ is
sufficiently large).

Now write
\[
r_1=\exp(a'/2)
\]
and then define
\[
r_n=\exp\( \frac{a'-(b_1+\cdots+b_{n-1})}{2} \)
\]
for all~$n\ge2$, where~$b_j$ is defined as in~\eqref{lisa}.

We wish to show by induction that, if~$X$ (and hence~$Y$) is chosen large
enough, then the following three properties hold for every~$n\ge1$:
\begin{equation}\label{andthenlaterwhenitgetsdark}
1<r_n< 1+\frac{\delta}{R_n},
\end{equation}
\begin{equation}\label{wegohome}
a'\(
1-\textstyle\frac{1}{2^n}-f(n)
\)
<
b_1+\cdots+b_n
<
a'\(
1-\textstyle\frac{1}{2^n}+f(n)
\)
\end{equation}
and
\begin{equation}\label{feedanimalsinthezoo}
\sum_{2^{R_n}<p\le 2^{R_nr_n}}\textstyle\frac{\log p}{p}\ll 2^{-n/2},
\end{equation}
where
\[
f(n)=\sum_{j=1}^{n}100^{-2^{j/4}}2^{j-n}.
\]
Notice that
\begin{equation}\label{thenlateramovietoo}
f(n)<2^{-(n+2)}
\end{equation}
for all~$n\ge1$, since
\[
\sum_{j=1}^{\infty}100^{-2^{j/4}}2^{j}<\textstyle\frac14.
\]
Thus~\eqref{uter},~\eqref{andthenlaterwhenitgetsdark},~\eqref{wegohome}
and~\eqref{feedanimalsinthezoo} together give the result.

If~$Y$ is large enough then, using~\eqref{uter}, we have
\[
1 < r_1=\exp(a'/2) < \(1+\frac{\delta}{Y}\)^{2/3}
< 1+\frac{7}{10}\cdot \frac{\delta}{Y}
< 1+\frac{\delta}{\left[2^{1/2}Y\right]} = 1+\frac{\delta}{R_1},
\]
and hence the bounds~\eqref{andthenlaterwhenitgetsdark} hold for~$n=1$.
Turning to~\eqref{wegohome}, notice that, since~$r_1=\exp(a'/2)$, we have
\[
b_1=\textstyle\frac{a'}{2}+\bigo\(\exp\(-C_{\ref{lastc}}2^{1/4}Y^{1/2}\)\)
\]
by~\eqref{lisa}, so~\eqref{wegohome} holds for~$n=1$ if~$Y$ is large
enough. Moreover,~\eqref{feedanimalsinthezoo} holds trivially for~$n=1$.

We assume now that~$X$ has been chosen large enough such that~\eqref{uter}
holds,~\eqref{wegohome} holds for~$n=1$, and
\begin{equation}\label{largeX}
 a' \exp\(C_{\ref{lastc}}Y^{1/2}\) > \frac{1}{5}\cdot \frac{\delta}{Y} \cdot
 \exp\(C_{\ref{lastc}}Y^{1/2}\)> 100 \quad \mbox{if } Y>X.
\end{equation}
In particular, the base step of the induction holds. Now assume
that~\eqref{andthenlaterwhenitgetsdark},~\eqref{wegohome}
and~\eqref{feedanimalsinthezoo} hold for some~$n=k-1$, with~$k\ge2$.
By~\eqref{wegohome} for~$n=k-1$, we have
\begin{equation}\label{problemsallleftalone}
a'\(\textstyle\frac{1}{2^k}-\frac{f(k-1)}{2}\)
< \log r_k
< a'\(\textstyle\frac{1}{2^k}+\frac{f(k-1)}{2}\).
\end{equation}
Using~\eqref{uter},~\eqref{thenlateramovietoo}
and~\eqref{problemsallleftalone}, we deduce that
\[
1 < r_k < \exp\(\frac{5}{4} \cdot \frac{a'}{2^k}\)
< \left(1+\frac{\delta}{Y}\right)^{(5/3) \cdot 2^{-k}}
< 1+\frac{7}{4} \cdot \frac{1}{2^k} \cdot \frac{\delta}{Y}
< 1+\frac{\delta}{\left[2^{k/2}Y\right]}
= 1+\frac{\delta}{R_k},
\]
and hence~\eqref{andthenlaterwhenitgetsdark} holds for~$n=k$.
Using~\eqref{lisa},~\eqref{thenlateramovietoo},~\eqref{problemsallleftalone}
and the definition of~$R_k$, we have
\begin{eqnarray*}
\sum_{2^{R_k}<p\le 2^{R_kr_k}}\frac{\log p}{p}
&\ll & \log 2^{R_kr_k}\sum_{2^{R_k}
\ <\ p\le 2^{R_kr_k}}\frac{1}{p} \\
&\ll & R_k\sum_{2^{R_k}<p\le 2^{R_kr_k}} \log\(1+\frac{1}{p}\)
\ \ll\  2^{-k/2}.
\end{eqnarray*}
It follows that~\eqref{feedanimalsinthezoo} holds for~$n=k$.
Moreover, by~\eqref{lisa} and the definition of~$r_k$, we have
\[
b_1+\cdots+b_k
= \frac{a'+b_1+\cdots+b_{k-1}}{2}
+ \bigo\(\exp\(-C_{\ref{lastc}}2^{k/4}Y^{1/2}\)\),
\]
and so
\begin{equation}\label{itssuchfun}
\frac{a'+b_1+\cdots+b_{k-1}}{2} - a'\cdot 100^{-2^{k/4}}
< b_1+\cdots+b_k
< \frac{a'+b_1+\cdots+b_{k-1}}{2} + a'\cdot 100^{-2^{k/4}},
\end{equation}
using~\eqref{largeX}.
From~\eqref{wegohome} for~$n=k-1$ and~\eqref{itssuchfun},
we deduce that
\[
a'\(1-\textstyle\frac{1}{2^k}-\frac{f(k-1)}{2}-100^{-2^{k/4}}\)
< b_1+\cdots+b_k
< a'\(1-\textstyle\frac{1}{2^k}+\frac{f(k-1)}{2}+100^{-2^{k/4}}\).
\]
This is equivalent to~\eqref{wegohome} for~$n=k$, since
\[
f(k)=\textstyle\frac{f(k-1)}{2}+100^{-2^{k/4}},
\]
completing the induction.
\end{proof}


\begin{thebibliography}{10}

\bibitem{MR0434929} Tom~M. Apostol.
\newblock {\em Introduction to analytic number theory}.
\newblock Springer-Verlag, New York, 1976.
\newblock Undergraduate Texts in Mathematics.

\bibitem{MR1461206} V.~Chothi, G.~Everest, and T.~Ward.
\newblock {$S$}-integer dynamical systems: periodic points.
\newblock {\em J. Reine Angew. Math.}, 489:99--132, 1997.

\bibitem{EMSW} G.~Everest, R.~Miles, S.~Stevens, and T.~Ward.
\newblock Orbit-counting in non-hyperbolic dynamical systems.
\newblock {\em J. Reine Angew. Math.}, 608:155--182, 2007.

\bibitem{EMSW2} G.~Everest, R.~Miles, S.~Stevens, and T.~Ward.
\newblock Dirichlet series for finite combinatorial rank dynamics.
\newblock {\em Trans. Amer. Math. Soc.}, 362(1):199--227, 2010.

\bibitem{MR2180241} G.~Everest, V.~Stangoe, and T.~Ward.
\newblock Orbit counting with an isometric direction.
\newblock In {\em Algebraic and topological dynamics}, volume 385 of {\em
  Contemp. Math.}, pages 293--302. Amer. Math. Soc., Providence, RI, 2005.

\bibitem{MR2135478} G.~Everest and T.~Ward.
\newblock {\em An introduction to number theory}, volume 232 of {\em Graduate
  Texts in Mathematics}.
\newblock Springer-Verlag London Ltd., London, 2005.

\bibitem{MR2000e:11087} G.~R. Everest and T.~Ward.
\newblock {\em Heights of polynomials and entropy in algebraic dynamics}.
\newblock Springer-Verlag London Ltd., London, 1999.

\bibitem{MR2647984} John Friedlander and Henryk Iwaniec.
\newblock {\em Opera de cribro}, volume~57 of {\em American Mathematical
  Society Colloquium Publications}.
\newblock American Mathematical Society, Providence, RI, 2010.

\bibitem{MR0186653} H.~Hasse.
\newblock \"{U}ber die {D}ichte der {P}rimzahlen {$p$}, f\"ur die eine
  vorgegebene ganzrationale {Z}ahl {$a\not=0$} von durch eine vorgegebene
  {P}rimzahl {$l\not=2$} teilbarer bzw. unteilbarer {O}rdnung {${\rm mod.}\,p$}
  ist.
\newblock {\em Math. Ann.}, 162:74--76, 1965/1966.

\bibitem{MR0205975} H.~Hasse.
\newblock \"{U}ber die {D}ichte der {P}rimzahlen {$p$}, f\"ur die eine
  vorgegebene ganzrationale {Z}ahl {$a\not=0$} von gerader bzw.ungerader
  {O}rdnung mod. {$p$}\ ist.
\newblock {\em Math. Ann.}, 166:19--23, 1966.

\bibitem{Sawian} S.~Jaidee.
\newblock {\em Mertens' theorem for arithmetical dynamical systems}.
\newblock PhD thesis, Univ. of East Anglia, 2010.

\bibitem{MR0068565d} Edmund Landau.
\newblock {\em Handbuch der {L}ehre von der {V}erteilung der {P}rimzahlen}.
\newblock Chelsea Publishing Co., New York, 1974.
\newblock 3rd ed, With an appendix by Paul T. Bateman.

\bibitem{MR0346130} D.~Lind.
\newblock Ergodic automorphisms of the infinite torus are {B}ernoulli.
\newblock {\em Israel J. Math.}, 17:162--168, 1974.

\bibitem{MR961739} D.~A. Lind and T.~Ward.
\newblock Automorphisms of solenoids and {$p$}-adic entropy.
\newblock {\em Ergodic Theory Dynam. Systems}, 8(3):411--419, 1988.

\bibitem{mertens} F.~Mertens.
\newblock Ein {B}eitrag zur analytyischen {Z}ahlentheorie.
\newblock {\em J. Reine Angew. Math.}, 78:46--62, 1874.

\bibitem{MR2378655} Hugh~L. Montgomery and Robert~C. Vaughan.
\newblock {\em Multiplicative number theory. {I}. {C}lassical theory},
  volume~97 of {\em Cambridge Studies in Advanced Mathematics}.
\newblock Cambridge University Press, Cambridge, 2007.

\bibitem{MR2486259} A.~Pakapongpun and T.~Ward.
\newblock Functorial orbit counting.
\newblock {\em J. Integer Seq.}, 12(2):Article 09.2.4, 20 pp., 2009.

\bibitem{MR727704} W.~Parry and M.~Pollicott.
\newblock An analogue of the prime number theorem for closed orbits of {A}xiom
  {A} flows.
\newblock {\em Ann. of Math. (2)}, 118(3):573--591, 1983.

\bibitem{MR1873399} Y.~Puri and T.~Ward.
\newblock Arithmetic and growth of periodic orbits.
\newblock {\em J. Integer Seq.}, 4(2):Article 01.2.1, 18 pp., 2001.

\bibitem{MR0143728} A.~Schinzel.
\newblock On primitive prime factors of {$a^{n}-b^{n}$}.
\newblock {\em Proc. Cambridge Philos. Soc.}, 58:555--562, 1962.

\bibitem{MR1139566} R.~Sharp.
\newblock An analogue of {M}ertens' theorem for closed orbits of {A}xiom {A}
  flows.
\newblock {\em Bol. Soc. Brasil. Mat. (N.S.)}, 21(2):205--229, 1991.

\bibitem{MR1342300} G.~Tenenbaum.
\newblock {\em Introduction to analytic and probabilistic number theory},
  volume~46 of {\em Cambridge Studies in Advanced Mathematics}.
\newblock Cambridge University Press, Cambridge, 1995.
\newblock Translated from the second French edition (1995) by C. B. Thomas.

\bibitem{MR2085157} T.~Ward.
\newblock Group automorphisms with few and with many periodic points.
\newblock {\em Proc. Amer. Math. Soc.}, 133(1):91--96, 2005.

\bibitem{MR736730} K.~Wiertelak.
\newblock On the density of some sets of primes. {IV}.
\newblock {\em Acta Arith.}, 43(2):177--190, 1984.

\bibitem{MR2422026} A.~J. Windsor.
\newblock Smoothness is not an obstruction to realizability.
\newblock {\em Ergodic Theory Dynam. Systems}, 28(3):1037--1041, 2008.

\end{thebibliography}

\end{document}